\input epsf
\magnification=\magstep1

\font\erm = cmr8

\font\smmmi=cmmi8
\font\smsy=cmsy8

\font\timmi=cmmi7

\def\small{\textfont1=\smmmi\textfont2=\smsy\textfont0=\erm\scriptfont1=\timmi}

\hsize=6.35truein \vsize=9truein

\font\title=cmssdc10 at 15pt

\def\frame #1{\vbox{\hrule height.1pt
\hbox{\vrule width.1pt\kern 10pt
\vbox{\kern 10pt
\vbox{\hsize 3.5in\noindent #1}
\kern 10pt}
\kern 10pt\vrule width.1pt}
\hrule height0pt depth.1pt}}
\def\bframe #1{\vbox{\hrule height.1pt
\hbox{\vrule width.1pt\kern 10pt
\vbox{\kern 10pt
\vbox{\hsize 5in\noindent #1}
\kern 10pt}
\kern 10pt\vrule width.1pt}
\hrule height0pt depth.1pt}}
\baselineskip 15pt

\centerline{\title Strength in numbers? Not always!}
\bigskip
\footnote{}{Key words and phrases: quorum sensing probability model}
\footnote{}{Partially supported by NSF grant DMS-0701396}

\bigskip
Rinaldo B. Schinazi

University of Colorado at Colorado Springs

email: rschinaz@uccs.edu
\bigskip
{\bf Abstract. } We propose a simple model to compute the probability of success under a quorum sensing strategy. We show that a quorum sensing strategy has a higher probability of success than an individualistic strategy when, for instance, the probability of success for a single individual is low and the cost of building a quorum is not too high.
On the other hand if the cost of building a quorum is too high then the probability of success under quorum sensing always decreases as a function of the quorum.  
\bigskip
{\bf 1. Introduction.} 'Quorum sensing' describes a strategy used by some bacteria under which the bacteria multiply until a critical mass is reached. At that point the bacteria turn on their virulence genes and launch an attack on their host. Several human diseases (such as cholera) are caused by quorum sensing bacteria. See Miller and Bassler (2001) and Zhu et al. (2002). All the literature we have found on the subject deals with the biology of quorum sensing and its evolutionary consequences. How do the bacteria know they have reached a quorum? Is this the pathway to multicellular organisms? See Winans and Bassler (2002).
What we have not seen addressed is why did quorum sensing appear? More precisely, is it obvious that there is always strength in numbers? We exhibit a simple probability model that shows that things may be less simple than they appear.
\bigskip
{\bf 2. The model.} We start with one individual. There are two phases in the process.
In the first phase either the population (of bacteria, for instance) disappears or it gets to the quorum $N$ where $N$ is a natural number. We will use a so called birth and death chain to model this first phase. If the population reaches $N$ then the second phase kicks in. In the second phase we assume that each one of the $N$ individuals has a probability $\rho$ of being successful. Success may mean different things in different situations. For a bacterium it may mean not being eliminated by the host. There is evidence that in the case of cholera not only do the bacteria multiply but they also evolve into more pathogenic strains before launching their attack on the host, see LaRocque et al. (2005). So we may think of $N$ as the number of strains rather than the number of individual bacteria and we may think of $\rho$ as the probability that a given strain escapes the immune system of the host. 
Let $A_N$ be the event that the population or the number of strains eventually reaches $N$ (starting with one). The alternative to the event $A_N$ is that the population gets killed off before reaching $N$.
In order for at least one individual in the population to be successful we first need the population to reach $N$ and then we need at least one of the $N$ individuals to be successful. The probability of the event just described is
$$f(N,\rho)=P(A_N)(1-(1-\rho)^N).$$
In order to get the expression above for $f$ we are assuming that the two phases are independent and that the $N$ individuals or $N$ strains present when the second phase kicks in are independently successful. The function $f$ represents the probability of success under a quorum sensing strategy, $N$ represents the quorum. We are interested in $f$ as a function of $N$. If there is always 'strength in numbers' then $f$ must be always increasing. In fact,  we will show that $f$ may be increasing, decreasing or neither.  

We use a birth and death chain to model the first phase. Assume that there are $n$ 
individuals at some point where $1\leq n<N$. Then there is a death with probability
$q_n$ or a birth with probability $p_n$. That is, for $n\geq 1$
$$\eqalign{n\to n-1 &\hbox{ with probability }q_n\cr
n\to n+1 &\hbox{ with probability }p_n\cr}$$
where $p_n+q_n=1$ for all $1\leq n<N$. 
We also assume that if the chain gets to 0 before getting to $N$ then it stays there.  There is a well known exact formula for the probability of eventually getting to $N$ when starting at 1. Namely, 
$$P(A_N)={1\over \sum_{i=0}^{N-1} P_i},$$
where $P_0=1$ and for $i\geq 1$
$$P_i=\Pi_{k=1}^i {q_k\over p_k}.$$
See, for instance, the formula at the bottom of page 16 in Schinazi (1999).
Next, we consider two particular birth and death chains. 
\medskip
{\bf 2.1 A monotone example.} In this subsection we assume that the birth and death probabilities are constant. That is, $p_n=p$ for all $n$ in $[1,N-1]$ and therefore $q_n=1-p=q$ for all such $n$. Let $r=q/p$. We get
$P_i=r^i,$
$$P(A_N)={r-1\over r^N -1}$$
and
$$f(N,\rho)=P(A_N)(1-(1-\rho)^N)={r-1\over r^N -1}(1-(1-\rho)^N).$$
A factorization yields
$$f(N,\rho)=\rho{1+(1-\rho)+\dots+(1-\rho)^{N-1}\over
1+r+\dots+r^{N-1}}.$$
Using this last expression and some simple algebra
it is not difficult to show that $f$ is decreasing as a function of $N$ if $r>1-\rho$ and increasing if $r<1-\rho$. Note that if $r\geq 1$ then $f$ is always decreasing. Hence, there is a dramatic change depending whether $\rho$ is smaller than or larger than $1-r$ provided $r<1$. This is illustrated by Figure 1 below for which we have $1-r$ a little above $0.18$. For $\rho=0.17$ we have an increasing function while for $\rho=0.19$ we have a decreasing one. However, note that in both cases the function $f$ becomes almost constant rather rapidly.
\bigskip
{\parindent = 0.5in
\tabskip=1em plus 2em minus 0.5em
\halign{\indent\hfil#\hfil && \hfil#\hfil \cr
{\epsfxsize=6cm\epsfbox{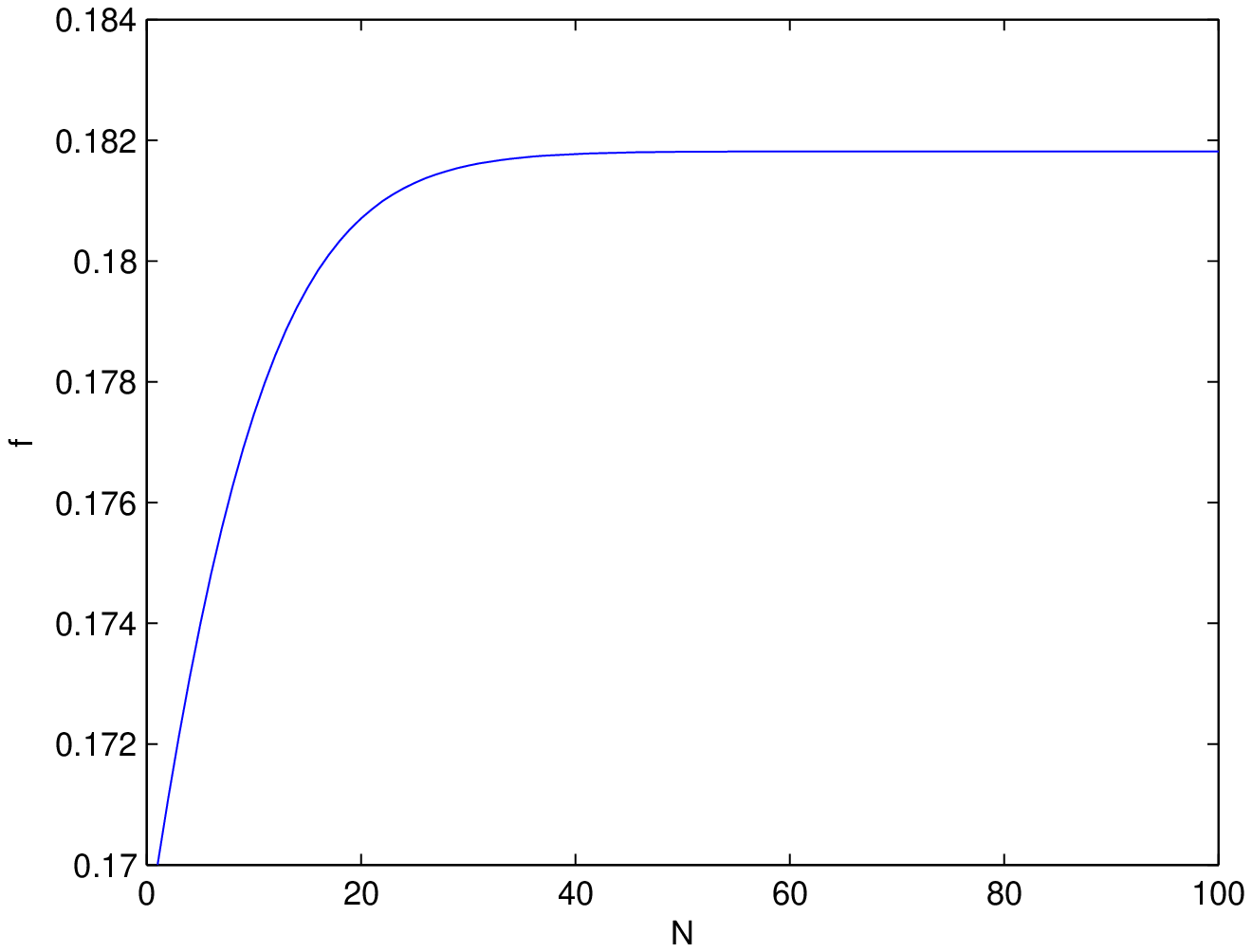}}& & {\epsfxsize=6cm\epsfbox{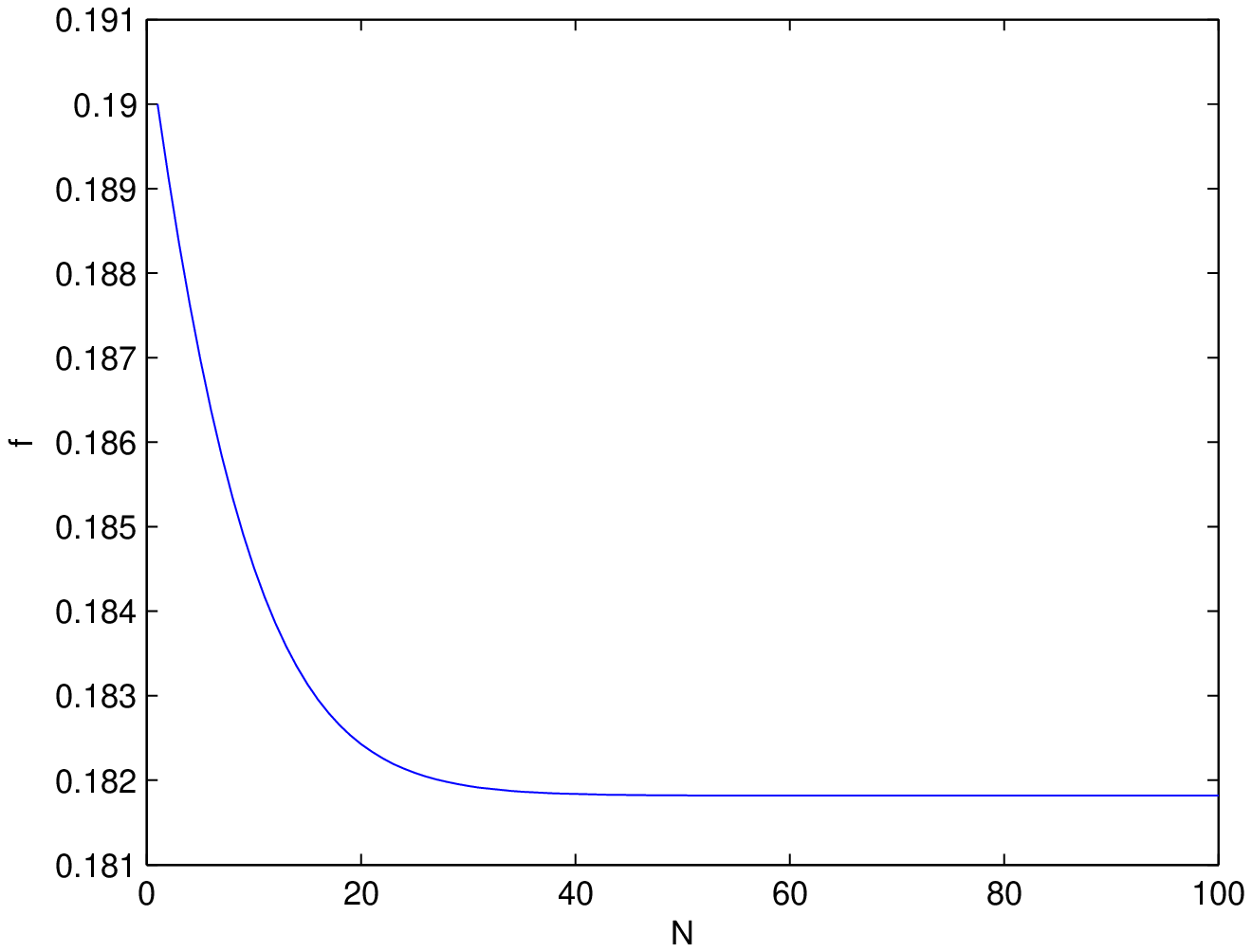}}\cr}}
\small{{\bf Figure 1.} These graphs of $f$ as a function of $N$ correspond to the birth and death chain of Section 2.1. We set $p=0.55$ and hence $q=0.45$. We picked $\rho=0.17$ for the graph on the left and $\rho=0.19$ for the graph on the right.}

\medskip
{\bf 2.2 A non monotone example.} In this subsection we consider the birth and death
chain with 
$$p_n={n+1\over 2n+1}\hbox{ and }q_n={n\over 2n+1},$$
for $n$ in $[1,N-1]$. In this case we get
$P_i={1\over i+1}$,
$$P(A_N)={1\over \sum_{i=1}^N {1\over i}},$$
and
$$f(N,\rho)={1\over \sum_{i=1}^N {1\over i}}(1-(1-\rho)^N).$$
For $\rho$ less than approximately 0.03 $f$ is increasing as a function of $N$. For $\rho$ larger than $0.5$ $f$ is decreasing. What makes this an interesting example is that for $\rho$ in $(0.03,0.5)$ the function $f$ is not monotone. See Figure 2 below. In this case it is for an intermediate value of $N$ that the probability of success under a quorum sensing strategy is maximal. 
 
\epsfxsize=6cm\epsfbox{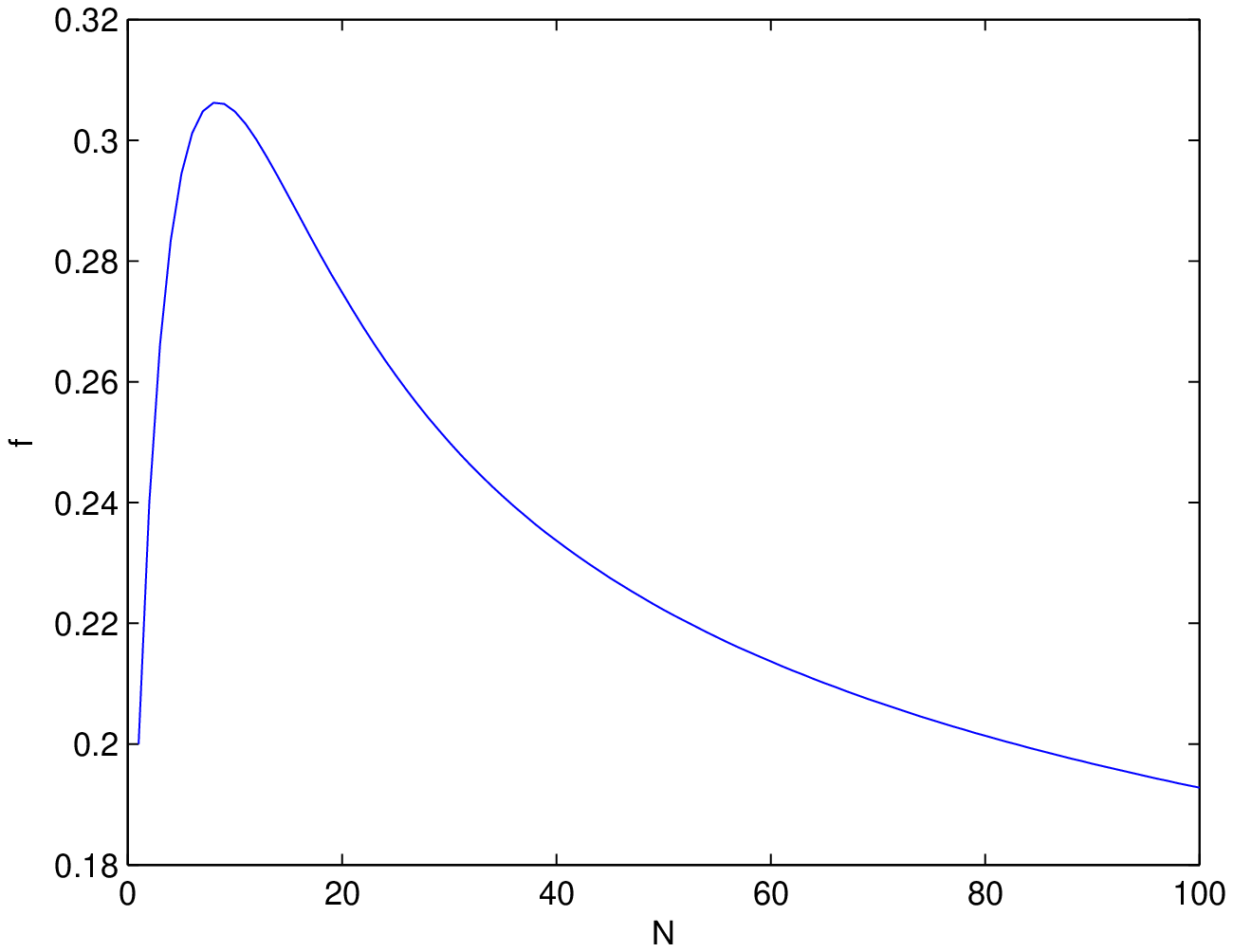}
\small{{\bf Figure 2.} This graph corresponds to the birth and death chain of
Section 2.2. We have set $\rho=0.2$.}
\bigskip
{\bf 3. Discussion.} Several types of bacteria use a quorum sensing strategy. Hence, this strategy must have some advantage over the 'every bacterium for itself' strategy in at least some situations. In this note we are interested in whether quorum sensing increases the probability of collective success. What we mean by 'collective success' is that at least one individual in the population be successful or in the case of bacteria that at least one of the strains escapes the immune system of the host. Maybe not surprisingly, quorum sensing is most effective when the probability of individual success $\rho$ is small. There are cases (see subsection 2.1 with $r<1$) for which there is a threshold for $\rho$ above which a quorum sensing strategy is less effective than an individualistic one. On the other hand if the cost of building a quorum is too high then the probability of success under quorum sensing decreases as a function of the quorum. This is, for instance, the case if $r\geq 1$ in the model of subsection 2.1. In that situation the probability of death is higher than the probability of birth and building a quorum is too costly. Hence, an individualistic strategy has a higher probability of success. 
Our model also shows that there are cases where quorum sensing is most effective for a an intermediate value of quorum (see subsection 2.2).  
\bigskip
{\bf References.}

R.C. LaRoque et al. (2005). Transcriptional profiling of Vibrio Cholerae recovered directly from patient specimens during early and late stages of human infection. Infection and Immunity 73, 4488-4493.

M.B. Miller and B.L. Bassler (2001). Quorum sensing in bacteria. Annu. Rev. Microbiol. 55, 165-199.

R.B. Schinazi (1999). Classical and spatial stochastic processes. Birkhauser.

S.C.Winans and B.L. Bassler (2002). Mob psychology. Journal of Bacteriology 184, 873-883.

J. Zhu, M.B. Miller, R.E. Vance, M. Dziejman, B.L. Bassler and J.J. Mekalanos (2002). PNAS 99, 3129-3134.

\end